\documentclass[twocolumn]{autart}    

\usepackage[utf8]{inputenc}
\usepackage[T1]{fontenc}
\usepackage[english]{babel}
\usepackage{graphicx}
\graphicspath{{Figures/}}

\usepackage[retainorgcmds]{IEEEtrantools}
\usepackage{mathtools}   
\usepackage{cases}
\usepackage{amssymb}
\usepackage{bbm}
\usepackage{enumerate}
\usepackage{siunitx}

\providecommand{\abs}[1]{\left\lvert#1\right\rvert}
\providecommand{\bigabs}[1]{\bigl\lvert#1\bigr\rvert}

\providecommand{\norm}[1]{\left\lVert#1\right\rVert}

\begin{document}

\begin{frontmatter}
\runtitle{Controllability of the 1D Schrödinger equation}  

\title{Controllability of the 1D Schrödinger equation using flatness} 


\author[First]{Philippe Martin}\ead{philippe.martin@mines-paristech.fr},  
\author[First,Second]{Lionel Rosier}\ead{lionel.rosier@mines-paristech.fr},
\author[First]{Pierre Rouchon}\ead{pierre.rouchon@mines-paristech.fr}

\address[First]{Centre Automatique et Systèmes, MINES ParisTech, PSL Research University\\
60 boulevard Saint-Michel, 75272 Paris Cedex 06, France}
\address[Second]{Centre de Robotique, MINES ParisTech, PSL Research University\\
	60 boulevard Saint-Michel, 75272 Paris Cedex 06, France}


\begin{abstract}                          
We derive in a direct way the exact controllability of the 1D Schrödinger equation with Dirichlet boundary control. We use the so-called flatness approach, which consists in parametrizing the solution and the control by the derivatives of a ``flat output''. This provides an explicit and very regular control achieving the exact controllability in the energy space.
\end{abstract}

\end{frontmatter}

\section{Introduction} The exact controllability of the linear Schrödinger equation (or of the related plate equation) is investigated in~\cite{lions-book,machtyngier,komornik-book} with the multiplier method, in~\cite{haraux,jaffard,KL-book} with nonharmonic Fourier analysis, in~\cite{lebeau} with microlocal analysis, and in~\cite{CFNS,liu} with frequency domain tests. It is extended to the semilinear Schrödinger equation in~\cite{RZ2009S1D,RZ2009SND,RZ2010} and~\cite{laurentS1D,laurentS3D}, by means of Strichartz estimates and Bourgain analysis. All those results are indirect, as they rely on observability inequalities for the adjoint system. A direct approach not involving the adjoint system is proposed in~\cite{LM,R2002,LT}; in particular, \cite{R2002} uses a fundamental solution of the Schrödinger equation with compact support in time, and provides controls with Gevrey regularity.

In this paper, we derive in a direct way the null (hence exact, since the equation is time-reversible) controllability of the Schrödinger equation with Dirichlet control,
\begin{IEEEeqnarray}{rCl'l}
	i\theta_t+\theta_{xx} &=& 0, &(t,x)\in (0,T)\times(0,1), \label{B1}\\
	\theta (t,0) &=& 0, & t\in (0,T), \label{B2}\\
	\theta (t,1) &=& u(t), & t\in (0,T), \label{B3}\\
	\theta (0,x) &=& \theta_0 (x), &x\in (0,1). \label{B4}
\end{IEEEeqnarray}
The initial condition~$\theta_0$, the state~$\theta$, and the control~$u$ are complex-valued functions. Given any final time $T>0$ and any initial condition~$\theta_0\in L^2(0,1)$, we explicitly build a regular control such that the state reached at~$T$ is exactly zero.

We use the so-called {\em flatness} approach, which consists in parametrizing the solution $\theta$ and the control~$u$
by the derivatives of a ``flat output'' $y$; this notion was initially introduced for finite-dimensional systems~\cite{FLMR},
and later extended to partial differential equations, see e.g.~\cite{LMR,PetitR2002TAC,PetitR2002SICON,LynchR2002IJC,WoittR2003COCV,MeurerThullKugi2008,Meurer2013}. A similar flatness-based approach is used in~\cite{MartiRR2014Automatica} for explicitly establishing the null controllability of the heat equation, and generalized to 1D parabolic equations in~\cite{MartiRR2016SICON}, with a control comprising two phases: a  first ``regularization'' phase, taking the irregular (namely, $L^2(0,1)$) initial state to a much more regular (namely, Gevrey of order~$\tfrac{1}{2}$) intermediate state; and a second phase using the flat parametrization to transfer this intermediate state to zero; a similar approach is used in~\cite{Moyan2016MCSS} in the case of a strongly degenerate parabolic equation. But whereas a zero control is sufficient in the first phase thanks to the natural smoothing effect of parabolic equations, the picture is very different for the Schrödinger equation, where no such smoothing effect exists with a zero control. Nevertheless, it is still possible to take the initial state to an intermediate state which is Gevrey of order~$\tfrac{1}{2}$, but with a well-chosen nonzero control; the idea, inspired by~\cite{R2002,LT,RZ2009SND} is to interpret the problem as a Cauchy problem on~$\Rset$, where some smoothing effect does take place; notice our result improves~\cite{LT}, with a more direct proof: the trajectory in the first phase is Gevrey of order~$1$ in time and~$\tfrac{1}{2}$ in space, compared to order~$2$ in time and~$1$ in space for~\cite{LT}. The second phase is then similar in spirit to the flatness-based control of~\cite{MartiRR2014Automatica,MartiRR2016SICON,Moyan2016MCSS}; with respect to~\cite{LT}, which solves an ill-posed problem in an abstract way, our control is explicitly given as a series.

Another contribution of this paper with respect to~\cite{MartiRR2014Automatica,MartiRR2016SICON} is a cleaner construction of the control, which emphasizes that the key point for establishing the null controllability of a flat partial differential equation is to find a trajectory which connects the initial condition to a sufficiently regular intermediate state; the control in the two phases is then automatically deduced from this trajectory. This construction also eliminates the constraint on the intermediate time that was needed in the preliminary version of this paper~\cite{MartiRR2014IFAC}.

The paper is organized as follows. In Section~\ref{sec:preliminaries}, we recall some preliminary notions. In Section~\ref{sec:flatness}, we precisely state what flatness means for the considered equation. In Section~\ref{sec:controllability}, we explicitly derive a control achieving null controllability.
In Section~\ref{sec:implementation}, we detail how this control can be numerically implemented. Finally, we illustrate in Section~\ref{sec:simulation} the effectiveness of the approach on a numerical example.

\section{Preliminaries}\label{sec:preliminaries}
We collect here a few definitions and facts that will be used throughout the paper.

We say  $y\in C^\infty([0,T])$ is {\em Gevrey of order~$s\geq0$ on~$[0,T]$} if there exist positive constants~$M,R$ such that
\begin{IEEEeqnarray*}{rCl}
\abs{y^{(p)}(t)} &\le& M\frac{p!^s}{R^p},~\forall t\in [0,T],\forall p\ge 0.
\end{IEEEeqnarray*}
More generally, if $K\subset\Rset^2$ is a compact set and $(t,x)\mapsto y(t,x)$ is a function of class~$C^\infty$ on~$K$ 
(i.e. $y$ is the restriction to $K$ of a function of class~$C^\infty$ on some open neighborhood $\Omega$ of~$K$), 
we say $y$ is {\em Gevrey of order~$(s_1,s_2)$ on~$K$} if there exist positive constants $M,R_1,R_2$ such that forall $(t,x)\in K$ and $(p_1,p_2)\in\Nset^2$,
\begin{IEEEeqnarray*}{rCl}
	\abs{\partial_{t}^{p_1}\partial_{x}^{p_2}y(t,x)} &\le&
	M\frac{p_1!^{s_1}p_2!^{s_2}}{R_1^{p_1}R_2^{p_2}}.
\end{IEEEeqnarray*}

By definition, a Gevrey function of order~$s$ is also of order~$r$ for~$r\geq s$. Gevrey functions of order~1 are analytic (entire if $s<1$).  Gevrey functions of order~$s>1$ may have a divergent Taylor expansion; the larger~$s$, the ``more divergent'' the Taylor expansion. Important properties of analytic functions generalize to Gevrey functions of order $s>1$: the scaling, addition, multiplication, inverse and derivation of Gevrey functions of order $s>1$ is of order~$s$, see e.g.~\cite{Yaman1989AGAG}. But contrary to analytic functions, Gevrey functions of order $s>1$ may be constant on an open set without being constant everywhere.  For example the ``step function''
\begin{IEEEeqnarray}{rCl}\label{eq:Gstep}
\phi_s(\rho)&:=&\begin{dcases*}
	1 & if $\rho\le0$\\
	0 & if $\rho\ge1$\\
	\dfrac{ e^{-\frac{M}{(1-\rho)^\sigma}} } { e^{-\frac{M}{\rho^\sigma}} + e^{-\frac{M}{(1-\rho)^\sigma}} } & if $\rho\in(0,1)$,
\end{dcases*}\end{IEEEeqnarray}
where $M>0$ and $\sigma:=(s-1)^{-1}$, is Gevrey of order~$s$ on~$\Rset$; notice 
$\phi_s^{(i)}(0)=\phi_s^{(i)}(1)=0$ for all~$i\geq1$. The result stems from the fact that $\rho\mapsto\mathbbm{1}_\Rset(\rho)e^{-\frac{1}{\rho^\sigma}}$ is Gevrey of order~$s$, which is classically proved thanks to the Cauchy integral formula, see \cite[chapter~$3.11$]{Widde1975book} for particular cases and \cite{LynchR2002IJC} for the general case.

From the two well-known properties of the Gamma function, see e.g.~\cite[section~$6.1$]{AbramS1992book},
\begin{IEEEeqnarray*}{rCl}
	\Gamma(\xi)\Gamma\Bigl(\xi+\frac{1}{2}\Bigr) &=& 2^{1-2\xi}\sqrt{\pi}~\Gamma(2\xi)\\
	\frac{\Gamma(\xi+a)}{\Gamma(\xi+b)} &\underset{\xi\to+\infty}{\sim}& \xi^{a-b},
\end{IEEEeqnarray*}
we deduce at once
\begin{IEEEeqnarray}{rCl}\label{eq:gamma}
	\Gamma(2\xi+1) &\underset{\xi\to+\infty}{\sim}& \frac{2^{2\xi}}{\sqrt{\pi\xi}}\Gamma^2(\xi+1).
\end{IEEEeqnarray}
Recall that $\Gamma(p+1)=p!$, when $p\in\Nset$. Following the usual notation in the literature on asymptotics,
\begin{IEEEeqnarray*}{rCl}
	f(\xi)\underset{\xi\to+\infty}{\sim}g(\xi) &\qquad\Leftrightarrow\qquad& \lim_{\xi\to+\infty}\frac{f(\xi)}{g(\xi)}=1.
\end{IEEEeqnarray*}

We will also use the obvious inequality
\begin{IEEEeqnarray}{rCl}\label{eq:fact}
	\frac{(p+q)!}{p!q!} &\le& 2^{p+q}.
\end{IEEEeqnarray}

\section{Flatness of the Schrödinger equation}\label{sec:flatness}

In this section we show the system~\eqref{B1}--\eqref{B4} is flat, with $y:=\theta_x(\cdot,0)$ as a flat output. Loosely speaking, it means there is a one-to-one correspondence between sufficiently regular solutions of ~\eqref{B1}--\eqref{B4} and sufficiently regular  arbitrary functions~$y$. This property will be paramount to prove controllability in the next section.

We first notice that, given any $C^\infty$ function~$Y$,  the formal series
\begin{IEEEeqnarray}{rCl}
	\theta (t,x) &=& \sum_{j\ge 0} \frac{x^{2j+1 } }{ (2j+1) ! } (-i)^j Y^{(j)} (t) \label{AA10a}\\
	u(t) &=&\sum _{j\ge 0} \frac{ (-i)^j}{(2j+1)!}Y^{(j)}(t) \label{AA10b}
\end{IEEEeqnarray}
obviously formally satisfy~\eqref{B1}--\eqref{B3} and 
\begin{IEEEeqnarray}{rCl}
	\theta_x(\cdot,0)=Y \label{eq:sol}.
\end{IEEEeqnarray}
The following property gives an actual meaning to these series when $Y$ is regular enough.

\begin{prop}\label{prop1}
	Let $s\in [0,2)$, $-\infty<t_1<t_2<\infty$, and assume $Y$ is Gevrey of order~$s$ on~$[t_1,t_2]$.
	Then $\theta$ defined by~\eqref{AA10a} is Gevrey of order~ $(s,\frac{s}{2})$ on $[t_1,t_2]\times[0,1]$; $u$ defined by~\eqref{AA10b} is Gevrey of order~$s$ on $[t_1,t_2]$.
\end{prop}

\begin{pf}
We want to prove the formal series
\begin{IEEEeqnarray}{rCl}\label{eq:dertheta}
 \partial_t^k\partial_x^l\theta(t,x) &=& \sum_{2j+1\ge l}\frac{x^{2j+1-l}}{(2j+1-l)!}(-i)^j Y^{(j+k)}(t)\IEEEeqnarraynumspace
\end{IEEEeqnarray}
is uniformly convergent on~$[t_1,t_2]\times [0,1]$,  with an estimate of its sum of the form
\begin{IEEEeqnarray*}{rCl}
 \abs{\partial_t^k\partial_x^l\theta(t,x)} &\le& C\frac{k!^s}{R_1^k}\, \frac{l!^\frac{s}{2}}{R_2^l}.
\end{IEEEeqnarray*}

By assumption, there are constants $M,R>0$ such that $\abs{Y^{(j)}(t)}\le M\frac{j!^s}{R^j}$ for all $j\geq0$ and all $t\in[t_1,t_2]$. Hence, for all $(t,x)\in[t_1,t_2]\times[0,1]$,
\begin{IEEEeqnarray*}{rCl}
	\abs{\frac{x^{2j+1-l}Y^{(j+k)}(t)}{(2j+1-l)!}} 
	&\le& \frac{M}{R^{j+k}}\, \frac{(j+k)!^s}{(2j+1-l)!}\\
	&\le& \frac{M}{R^{j+k}}\, \frac{(2^{j+k}j!k!)^s}{(2j+1-l)!}\\
	&\le& \frac{M'k!^s}{R_1^{j+k}}\, \frac{\bigl(2^{-2j}\sqrt{\pi(j+1)}\, (2j)!\bigr)^\frac{s}{2}}{(2j+1-l)!}\\
	&\le& \frac{M''k!^s}{R_1^{j+k}}\, \frac{\bigl(2^{-2j}\sqrt{\pi(j+1)}\, (2j+1)!\bigr)^\frac{s}{2}}{(2j+2)^\frac{s}{2}(2j+1-l)!}\\
	&\le& \frac{M''k!^s}{R_1^{j+k}}\, \frac{\bigl(\sqrt{\pi} \, (2j+1-l)!\,l!\bigr)^\frac{s}{2}}{(j+1)^\frac{s}{4}(2j+1-l)!}\\
	&=& M''A_{j,l}l!^{\frac{s}{2}}\, \frac{k!^s}{R_1^k},
\end{IEEEeqnarray*}
for some constants $M''\ge M'\ge M$; we have set $R_1:=2^{-s}R$ and
\begin{IEEEeqnarray*}{rCl}
	A_{j,l} &:=& \frac{\pi^{\frac{s}{4}}}{R_1^j(j+1)^\frac{s}{4}(2j+1-l)!^{1-\frac{s}{2}}},
\end{IEEEeqnarray*}
used \eqref{eq:fact} twice and \eqref{eq:gamma} for $\xi:=j$.
As $\sum_{2j+1\ge l}A_{j,l}<\infty$ (ratio test), the series in \eqref{eq:dertheta} is uniformly convergent for all $k,l\ge0$, implying $\theta\in C^\infty([t_1,t_2]\times[0,1])$. Moreover,
\begin{IEEEeqnarray*}{rCl}
	\sum_{2j+1\ge l} A_{j,l} 
	&=& \frac{(2\pi)^\frac{s}{4}\sqrt{R_1}}{\sqrt{R_1}^l}\sum_{p\ge0}\frac{1}{\sqrt{R_1}^p(p+l+1)^\frac{s}{4}p!^{1-\frac{s}{2}}}\\
	&\le& C(R_1)\frac{1}{\sqrt{R_1}^l}
\end{IEEEeqnarray*}
for $C(R_1)$ large enough. Hence,
\begin{IEEEeqnarray*}{rCl}
 \abs{\partial_t^k\partial_x^l\theta(t,x)} &\le& M''C(R_1)\,\frac{k!^s}{R_1^k}\, \frac{l!^\frac{s}{2}}{\sqrt{R_1}^l}.
\end{IEEEeqnarray*}
which proves $\theta$ is Gevrey of order~$(s,\frac{s}{2})$ on $[t_1,t_2]\times[0,1]$. As a simple consequence, $u=\theta(\cdot,1)$ is Gevrey of order~$s$ on $[t_1,t_2]$.\qed
\end{pf}

We can now give a precise meaning to \eqref{B1}--\eqref{B4} being flat with $\theta_x(\cdot,0)$ as a flat output: 
\begin{itemize}
	\item obviously, any solution $\theta$ of~\eqref{B1}-\eqref{B2} that is Gevrey of order~$(s,\frac{s}{2})$ on $[0,T]\times[0,1]$ --which implies $u$ is Gevrey of order~$s$ on $[0,T]$ and $\theta_0$ is Gevrey of order~$\frac{s}{2}$ on~$[0,1]$-- uniquely determines the flat output $\theta_x(\cdot,0)$ as a function that is Gevrey of order~$s$ on~$[0,T]$
	\item conversely, any function $Y$ Gevrey or order~$s$ on~$[0,T]$ gives rise by~\eqref{AA10a} to a solution $\theta$ of~\eqref{B1}-\eqref{B2}-\eqref{eq:sol}, which is Gevrey or order~$(s,\frac{s}{2})$ on $[0,T]\times[0,1]$; moreover, this solution is unique by the Holmgren theorem.
\end{itemize}
\section{Controllability of the Schrödinger equation}\label{sec:controllability}
We derive in the following Theorem~\ref{thm:Main} an explicit control steering the system from any initial state~$\theta_0\in L^2(0,1)$ at time~$0$ to the final state~$0$ at time~$T$. This proves null controllability, hence exact controllability since the Schrödinger equation is reversible with respect to time. This is a two-step procedure:
\begin{itemize}
	\item we first apply a ``smoothing'' control till some arbitrary intermediate time~$\tau<T$ so as to reach a ``smooth'' intermediate state~$\theta(\tau,\cdot)$
	\item we then apply a flatness-based control using the parametrization~\eqref{AA10b}.
\end{itemize}

Notice that the initial-boundary value problem for the Schrödinger equation~\eqref{B1}--\eqref{B4} with zero control does not smooth an $L^2$ initial condition into a state that is Gevrey of order~$\frac{s}{2}$ (all its eigenvalues are on the imaginary axis), which is in sharp contrast to the heat equation, see~\cite{MartiRR2014Automatica}. The control in the first phase, which therefore cannot be zero, is here obtained by using the dispersive effect of the Schrödinger equation on~$\Rset$, according to the following proposition.
\begin{prop}\label{prop2}
Let $\theta_0\in L^2(0,1)$, and define its odd extension $\theta_0^{odd}\in L^1(\Rset)\cap L^2(\Rset)$ by
\begin{IEEEeqnarray*}{rCl}
	\theta_0^{odd}(x) &:=& 
	\begin{cases*}
		\theta_0(x) & if $x\in(0,1)$\\		
		-\theta_0(-x) & if $x\in(-1,0)$\\
		0 & if $x\in(-\infty,-1]\cup[1,+\infty)$.
	\end{cases*}
\end{IEEEeqnarray*}
Introduce also the fondamental solution of the Schrödinger equation $E(t,x):=\displaystyle\frac{e^{i\frac{x^2}{4t}}}{\sqrt{4\pi it}}$.
Then the function defined by $\theta^-(0,\cdot):=\theta^{odd}_0$ and
\begin{IEEEeqnarray*}{rCl'l}
	\theta^-(t,x) &:=& \int_{-1}^1E(t,x-y)\theta^{odd}_0(y)dy, &(t,x)\in\Rset^*\times\Rset
\end{IEEEeqnarray*}
enjoys the following properties:
\begin{enumerate}[(i)]
	\item\label{it:c} $\theta^-$ belongs to $C^\infty(\Rset^*\times\Rset)$
	\item\label{it:f} $\theta^-$ is Gevrey of order~$(1,\frac{1}{2})$  on $[t_1,t_2]\times[-L,L]$ for all $0<t_1<t_2$ and~$L>0$.
	\item\label{it:d} $\theta^-(t,-x)=-\theta^-(t,x)$ for all $(t,x)\in\Rset^*\times\Rset$
	\item\label{it:a} $\theta^-|_{(0,T)\times(0,1)}$ is the solution of~\eqref{B1}--\eqref{B4}, where the control in~\eqref{B3} is given by $u(t):=\theta^-(t,1)$
	\item\label{it:b} $t\mapsto\theta^-(t,\cdot)$ belongs to~$C\bigl(\Rset,L^2(\Rset)\bigr)$
	\item\label{it:g} $\norm{\theta^-(t,\cdot)}_{L^2(\Rset)}=\sqrt2\norm{\theta_0}_{L^2(0,1)}$ for all $t\in\Rset$
	\item\label{it:e} $t\mapsto\norm{\theta^-(t,\cdot)}_{L^\infty(\Rset)}$ belongs to~$L^4(\Rset)$.
\end{enumerate}
\end{prop}
\begin{pf}
	We start by estimating the growth of the derivatives of~$E$; notice $E_t=iE_{xx}$ (since $E$ is the fundamental solution of~\eqref{B1}).
	It is well-known, see e.g.~\cite[$(5.5.3)$]{Szego1975book}, that
	\begin{IEEEeqnarray*}{rCl}
		\frac{d^k}{dz^k}e^{-z^2} &=& (-1)^ke^{-z^2}H_k(z),\label{eq:Hgenerating}
	\end{IEEEeqnarray*}
	where the $H_k$ are the (physicists') Hermite polynomials, which satisfy  $H_0(x)=1$, $H_1(x)=2x$ and
	\begin{IEEEeqnarray*}{rCl}
		H_{k+1}(x) &=& 2xH_{k}(x)-2kH_{k-1}(x).\label{eq:Hrecurrence}
	\end{IEEEeqnarray*}
	As a straightforward consequence,
	\begin{IEEEeqnarray}{rCl}\label{eq:Ek}
	\partial_x^kE(t,x) &=& \frac{(-1)^k}{(4it)^\frac{k}{2}}H_k\biggl(\frac{x}{\sqrt{4it}}\biggr)E(t,x)\nonumber\\
	&=& \frac{(-1)^k}{\sqrt{\pi}(4it)^\frac{k+1}{2}}e^{-\bigl(\frac{x}{\sqrt{4it}}\bigr)^2}H_k\biggl(\frac{x}{\sqrt{4it}}\biggr).
\end{IEEEeqnarray}
	On the other hand, see e.g.~\cite[Theorem~$8.22.7$]{Szego1975book},
\begin{IEEEeqnarray*}{rCl}
	\frac{\Gamma(\frac{k}{2}+1)}{\Gamma(k+1)}e^{-\frac{z^2}{2}}H_k(z) &\underset{k\to\infty}{\sim}& \cos\Bigl(z\sqrt{2k+1}-\frac{k\pi}{2}\Bigr)
\end{IEEEeqnarray*}
	uniformly for $\abs{z}\leq A$ whatever the positive constant~$A$. As an immediate consequence,
	\begin{IEEEeqnarray}{rCl}\label{eq:Hk}
	\abs{e^{-\frac{z^2}{2}}H_k(z)} &\le& M(A)\frac{\Gamma(k+1)}{\Gamma(\frac{k}{2}+1)}e^{A\sqrt{2k}}
\end{IEEEeqnarray}
	for $M(A)$ large enough. Taking $z:=\frac{x}{\sqrt{4it}}$ with $\abs{x}\le L$ and $\abs{t}\ge t_1$, hence $\abs{z}\le A:=\frac{L}{\sqrt{4\pi t_1}}$, we then get
	\begin{IEEEeqnarray*}{rCl}
	\abs{\partial_x^kE(t,x)} &=& \frac{ \bigabs{e^{-\frac{z^2}{2}}} }{\sqrt{\pi}\abs{4t}^\frac{k+1}{2}}\abs{e^{-\frac{z^2}{2}}H_k(z)}\\
	&\le& \frac{ e^{\frac{A^2}{2}} }{\sqrt{\pi}\abs{4t_1}^\frac{k+1}{2}}M(A)\frac{\Gamma(k+1)}{\Gamma(\frac{k}{2}+1)}e^{A\sqrt{2k}}\\
	&\underset{k\to\infty}{\sim}& \frac{ e^{\frac{A^2}{2}}M(A) }{(2\pi)^\frac{3}{4}\sqrt{t_1}}\frac{k!^\frac{1}{2}}{\sqrt{2t_1}^k}
	\frac{ e^{A\sqrt{2k}} }{k^\frac{1}{4}}\\
	&\le& M'(t_1,A)\frac{k!^\frac{1}{2}}{R^k},
\end{IEEEeqnarray*}
	for any $R\in(0,\sqrt{2t_1})$ and $M'(t_1,A)$ large enough; the third line comes from~\eqref{eq:gamma} with $\xi:=\frac{k}{2}$.
	As $E_t=iE_{xx}$,
\begin{IEEEeqnarray}{rCl}
	\abs{\partial_x^k\partial_t^lE(t,x)} &=& \abs{\partial_x^{k+2l}E(t,x)}\nonumber\\
	&\le& M'(t_1,A)\frac{(k+2l)!^\frac{1}{2}}{R^{k+2l}}\nonumber\\
	&\le& M'(t_1,A)\frac{2^{\frac{k}{2}+l}k!^\frac{1}{2}(2l)!^\frac{1}{2}}{R^{k+2l}}\nonumber\\
	&\le& M''(t_1,A)\frac{k!^\frac{1}{2}}{\bigl(\frac{R}{\sqrt2}\bigr)^k}\frac{l!}{\bigl(\frac{R}{2}\bigr)^{2l}},\label{eq:majE}
\end{IEEEeqnarray}
for $M''(t_1,A)$ large enough; we have used~\eqref{eq:fact} to get the second line, and \eqref{eq:gamma} with $\xi:=l$ to get the last line.

	This implies that for all compact set $K$ of~$\Rset^*\times\Rset$ and $k,l\in\Nset$, there is a constant $M_{K,k,l}>0$ such that
	\begin{IEEEeqnarray*}{rCl}
		\sup_{\substack{(t,x)\in K\\ y\in[-1,1]}}\abs{\partial_x^k\partial_t^lE(t,x-y)} &\le& M_{K,k,l}.
	\end{IEEEeqnarray*}
	We can now establish~\eqref{it:c} by a standard argument of differentiation under the integral sign. Indeed:
	\begin{itemize}
		\item the function $y\mapsto E(t,x-y)\theta^{odd}_0(y)$ belongs to $L^1(-1,1)$ for all $(t,x)\in\Rset^*\times\Rset$
		\item the function $(t,x)\mapsto E(t,x-y)\theta^{odd}_0(y)$ belongs to $C^\infty(\Rset^*\times\Rset)$ for almost all $y\in[-1,1]$
		\item $\abs{\partial_x^k\partial_t^lE(t,x-y)\theta^{odd}_0(y)} \le M_{K,k,l}\abs{\theta^{odd}_0(y)}$ for almost all $y\in[-1,1]$ and  all $(t,x)\in K$.
	\end{itemize}
	
	We then have~\eqref{it:f}, since by~\eqref{eq:majE},
	\begin{IEEEeqnarray*}{rCl}
		\abs{\partial_x^k\partial_t^l\theta^-(t,x)} &=& \abs{\int_{-1}^1\partial_x^{k+2l}E(t,x-y)\theta^{odd}_0(y)dy}\\
		&\le& 2M'(t_1,L+1)\norm{\theta_0}_{L^1(0,1)}\frac{k!^\frac{1}{2}}{\bigl(\frac{R}{\sqrt2}\bigr)^k}\frac{l!}{\bigl(\frac{R}{2}\bigr)^{2l}}
	\end{IEEEeqnarray*}
	for all $(t,x)\in[t_1,t_2]\times[-L,L]$.
	
	Statement~\eqref{it:d} follows from the change of variable $z:=y+x$ applied to~$\theta^-(t,-x)$ and the definition of~$\theta_0^{odd}$.
	
	$\theta^-$ is the (mild) solution of the Cauchy problem
	\begin{IEEEeqnarray*}{rCl'l}
		i\theta_t+\theta_{xx} &=& 0, &(t,x)\in\Rset^2\\
		\theta(0,x) &=& \theta_0^{odd}(x), & x\in\Rset,
	\end{IEEEeqnarray*}
	see~\cite[section~$4.1$]{LP-book} or \cite[section~$2.1$]{Cazenave2003book}; in particular \eqref{it:b} holds true. Obviously $\theta^-(t,0)=0$ and $\theta^-(t,1)=u(t)$ for~$t>0$, which implies~\eqref{it:a}.
	
	For all $t\in\Rset$, $\norm{\theta^-(t,\cdot)}_{L^2(\Rset)}=\norm{\theta_0^{odd}}_{L^2(\Rset)}$, see~\cite[section~$4.1$]{LP-book} or \cite[section~$2.1$]{Cazenave2003book}; \eqref{it:g} then obviously follows from the definition of~$\theta_0^{odd}$.
	
	Finally, \eqref{it:e} is Strichartz's estimate for the admissible pair $(4,+\infty)$, see~\cite[section~$4.2$]{LP-book} or \cite[section~$2.3$]{Cazenave2003book}.
\qed\end{pf}

\begin{thm}\label{thm:Main}
Let $\theta_0\in L^2(0,1)$ and $T>0$ be given. Pick any $\tau\in(0,T)$ and any~$s\in(1,2)$. Define
\begin{IEEEeqnarray*}{rCl}
	Y(t) &:=& \phi_s\Bigl(\frac{t-\tau}{T-\tau}\Bigr)\theta^-_x(t,0),
\end{IEEEeqnarray*}
with $\theta^-$ as in Proposition~\ref{prop2} and $\phi_s$ a Gevrey step function of order~$s$ such as~\eqref{eq:Gstep}.
Then the control 	
\begin{numcases}{u(t):=}
		\theta^-(t,1) & if $0\le t\le\tau$,\label{eq:control1}\\
		\sum_{j\ge0}\frac{(-i)^j}{(2j+1)!}Y^{(j)}(t) & if $\tau <t\le T$\qquad\label{eq:control2}
\end{numcases}
steers~\eqref{B1}--\eqref{B4} from~$\theta_0$ to~$\theta(T,\cdot)=0$; $u$ belongs to~$L^4(0,T)$, and is Gevrey of order~$s$ on $[\varepsilon,T]$ for all $\varepsilon\in(0,T)$.
The corresponding solution of~\eqref{B1}--\eqref{B4} is given by
\begin{numcases}{\theta(t,x):=}
\theta^-(t,x) & if $0\le t\le\tau$,\label{eq:traj1}\\
\sum_{j\ge0}\frac{(-i)^jx^{2j+1}}{(2j+1)!}Y^{(j)}(t) & if $\tau<t\le T$;\qquad\label{eq:traj2}
\end{numcases}
$t\mapsto\theta(t,\cdot)$ belongs to~$C\bigl([0,T],L^2(0,1)\bigr)$, and $\theta$ is Gevrey of order~$(s,\frac{s}{2})$ on
$[\varepsilon,T]\times[0,1]$ for all $\varepsilon\in(0,T)$.

\end{thm}
\begin{rem}
Notice the flatness property is used twice in the design of the control: of course in the very definition of the second phase through the parametrization~\eqref{AA10b}; and also in the definition of $Y(t)$ through $\theta^-_x(\cdot,0)$, which is the value of the flat output as computed in the first phase. This shows that the key point for establishing the null controllability of a flat partial differential equation is to find a trajectory which connects the initial condition to a sufficiently regular intermediate state; the control in the two phases is then automatically deduced from this trajectory. Though this is not completely apparent, this argument underlies the construction of the control in~\cite{MartiRR2014Automatica,MartiRR2016SICON}. Notice also that the construction of Theorem~\ref{thm:Main} is much simpler than the construction in the preliminary version~\cite{MartiRR2014IFAC}, as it avoids the additional series in Lemma~$5$, and removes the restriction $\tau>\tfrac{2}{3}T$.
\end{rem}
\begin{pf}
By Proposition~\ref{prop2}(\ref{it:a}), $\theta^-$ is the solution of~\eqref{B1}--\eqref{B4} on~$[0,\tau]\times[0,1]$; by~(\ref{it:b}), $t\mapsto\theta^-(t,\cdot)$ is in~$C\bigl(\Rset,L^2(\Rset)\bigr)$, and a fortiori in~$C\bigl([0,\tau],L^2(0,1)\bigr)$.  It eventually belongs to $C\bigl([0,T],L^2(0,1)\bigr)$ as an obvious consequence of its Gevrey regularity away from~$0$, to be proved in the sequel.
	
By Proposition~\ref{prop2}(\ref{it:f}), $\theta^-$ is Gevrey of order~ $(1,\frac{1}{2})$ on $[\varepsilon,T]\times[0,1]$ for all $\varepsilon\in(0,T)$, and a fortiori Gevrey of order~$(s,\frac{s}{2})$ since~$s>1$. In particular, $\theta^-(t,\cdot)$ is Gevrey of order~$\frac{1}{2}$, hence entire, for $t\in[\varepsilon,T]$; since it is moreover odd by~Proposition~\ref{prop2}(\ref{it:d}), it can be expanded as
\begin{IEEEeqnarray}{rCl}
	\theta^-(t,x) 
	&=& \sum_{j\ge0}\frac{x^{2j+1}}{(2j+1)!}\partial_x^{2j+1}\theta^-(t,0)\nonumber\\
	&=& \sum_{j\ge0}\frac{(-i)^jx^{2j+1}}{(2j+1)!}\partial_t^j\theta^-_x(t,0).\label{eq:entire}
\end{IEEEeqnarray}
Also $\theta^-_x(\cdot,0)$ is Gevrey of order~$s$ on $[\varepsilon,T]$; and so is~$Y$, as the product of two Gevrey functions of order~$s$. By Proposition~\ref{prop1}, the series in~\eqref{eq:traj2} is Gevrey of order~$(s,\frac{s}{2})$ on $[\varepsilon,T]\times[0,1]$ and satisfy~\eqref{B1}-\eqref{B3}. Moreover, \eqref{eq:traj1} and \eqref{eq:traj2} agree on~$[\varepsilon,\tau]\times[0,1]$. Indeed, $\phi_s\bigl(\frac{t-\tau}{T-\tau}\bigr)=1$ for $t\le\tau$, hence $Y(t)=\theta^-_x(t,0)$; by~\eqref{eq:entire}, this implies
\begin{IEEEeqnarray*}{rCl}
	\sum_{j\ge0}\frac{(-i)^jx^{2j+1}}{(2j+1)!}Y^{(j)}(t) &=& \theta^-(t,x).
\end{IEEEeqnarray*}
Therefore, $\theta$ given by~\eqref{eq:traj1}-\eqref{eq:traj2} is Gevrey of order~$(s,\frac{s}{2})$ on
$[\varepsilon,T]\times[0,1]$, and is the solution of~\eqref{B1}--\eqref{B4}. Finally, $\theta(T,\cdot)=0$; indeed, $\phi_s^{(m)}(1)=0$ for~$m\in\Nset$, which implies $Y^{(j)}(T)=0$ for~$j\in\Nset$.

Since the control~$u$ defined by~\eqref{eq:control1}-\eqref{eq:control2} is equal to $\theta(\cdot,1)$, it is Gevrey of order~$s$ on~$[\varepsilon,T]$. Moreover, Proposition~\ref{prop2}(\ref{it:e}) implies in particular that $u$ is in $L^4[0,\tau]$, hence in~$L^4[0,T]$.
\qed\end{pf}

\section{Numerical implementation of the control}\label{sec:implementation}
Theorem~\ref{thm:Main} provides an explicit expression for the control. We now detail how this expression can be effectively implemented in practice. The fact that a small error on the so computed control yields a small error on the final state follows from classical semigroup arguments.

\subsection{First phase}
The control in the first phase is given by the integral
\begin{IEEEeqnarray}{rCl}
	\theta^-(t,1) &:=& \frac{1}{\sqrt{4\pi it}}\int_{-1}^1 e^{i\frac{(y-1)^2}{4t}}\theta^{odd}_0(y)dy.\label{eq:uFirst}
\end{IEEEeqnarray}
This integral is difficult to evaluate numerically when $t$ gets small, because of the oscillating singularity at~$t=0$. A much better alternative is to expand it asymptotically.
We thus recall the basics of asymptotic expansion of integrals of the form
\begin{IEEEeqnarray*}{rCl's}
	I(\lambda) &:=& \int_a^bf(y)e^{i\lambda(y-\bar y)^2}dy, &$a<b\le\bar y$,
\end{IEEEeqnarray*}
when the real parameter $\lambda$ tends to $+\infty$, and refer the reader to~\cite[chapter~6]{BleistH1986book} for a detailed account.
The fundamental result is the following: consider the integral
\begin{IEEEeqnarray*}{rCl}
	J(\lambda) &:=& \int_0^{+\infty}F(u)e^{\varepsilon i\lambda u}du,
\end{IEEEeqnarray*}
where $\varepsilon=\pm1$ and $F$ is such that
\begin{itemize}
	\item $F(u)=0$ when $u>\bar u$, for some $\bar u>0$
	\item $F\in C^M(0,+\infty)$ for some $M\in\Nset$
	\item as $u\to0^+$, $F$ has the asymptotic expansion
	\begin{IEEEeqnarray}{rCl}
		F(u) &=& \sum_{m=0}^Mp_mu^{a_m}+o(u^{a_M}),\label{eq:Fasymptotic}
	\end{IEEEeqnarray}
	with $\Re(a_0)<\cdots<\Re(a_M)$; as usual, the ``little-o'' symbol in $\xi=o(\zeta)$ is used to mean $\frac{\xi}{\zeta}\to0$
	\item the asymptotic expansion as $u\to0^+$ of $F^{(1)},\ldots,F^{(M)}$ is obtained by differentiating~\eqref{eq:Fasymptotic} term-by-term.
\end{itemize}
Then, as $\lambda\to+\infty$, $J(\lambda)$ has the asymptotic expansion
\begin{IEEEeqnarray}{rCl}
	J(\lambda) &=& \sum_{m=0}^M\frac{p_m\Gamma(1+a_m)}{\lambda^{1+a_m}}e^{\varepsilon i\frac{\pi}{2}(1+a_m)}+o(\lambda^{-1-a_M}).\IEEEeqnarraynumspace\label{eq:fundamental}
\end{IEEEeqnarray}
This result can be formally obtained by replacing $F$ by its asymptotic expansion~\eqref{eq:Fasymptotic} in the integral and integrating term-by-term; the difficult part is to justify this process, since \eqref{eq:Fasymptotic} makes sense only near~$0$ and certainly not in the whole range of integration.

To evaluate $I(\lambda)$, we first ``isolate'' the endpoints of integration, see~\cite[section~3.3]{BleistH1986book}: we take a function~$\chi_a\in C^\infty(\Rset)$ which is identically~1 around~$a$ and identically~0 around~$b$ (such a function is for instance readily built from the Gevrey step function~$\phi_s$~\eqref{eq:Gstep}), and set $\chi_b:=1-\chi_a$. We then have $I(\lambda)=I_a(\lambda)+I_b(\lambda)$, with
\begin{IEEEeqnarray*}{rCl's}
	I_a(\lambda) &:=& \int_a^b\chi_a(y)f(y)e^{i\lambda(y-\bar y)^2}dy\\
	I_b(\lambda) &:=& \int_a^b\chi_b(y)f(y)e^{i\lambda(y-\bar y)^2}dy.
\end{IEEEeqnarray*}

By suitable changes of variables, $I_a(\lambda)$ and $I_b(\lambda)$ are next reduced to forms similar to $J(\lambda)$. 
Setting $r:=b-y$ and $u:=r(2\bar b+r)$, with $\bar b:=\bar y-b$, readily yields
\begin{IEEEeqnarray*}{rCl's}
	I_b(\lambda) &=& e^{i\lambda\bar b^2}\int_0^{+\infty}\frac{(\chi_bf)(b-r)}{2(\bar b+r)}\biggr|_{r=-\bar b+\sqrt{\bar b^2+u}}e^{i\lambda u}du.
\end{IEEEeqnarray*}
Taking the asymptotic expansion of the integrand as $u\to0^+$, we can now use the fundamental result~\eqref{eq:fundamental}. Finally, since $\chi_bf=f$ around~$b$, $(\chi_bf)(b-r)$ and~$f(b-r)$ have the same asymptotic expansion as $r\to0^+$, which means the asymptotic expansion of $I_b(\lambda)$ as $\lambda\to+\infty$ does not depend on $\chi_bf$, but only on~$f$. The process is similar with $I_a$: we set $s:=t-a$ and $u:=s(2\bar a-s)$, with $\bar a:=\bar y-a$, which yields
\begin{IEEEeqnarray*}{rCl's}
	I_a(\lambda) &=& e^{i\lambda\bar a^2}\int_0^{+\infty}\frac{(\chi_af)(a+s)}{2(\bar a-s)}\biggr|_{s=\bar a-\sqrt{\bar a^2-u}}e^{i\lambda u}du,
\end{IEEEeqnarray*}
and then use~\eqref{eq:fundamental}.

From the previous discussion, $I(\lambda)$ can be asymptotically expanded if
$f$ is smooth enough on $(a,b)$ and has well-behaved asymptotic expansions (similar to those required for $F$) as $y\to a^+$ and $y\to b^-$. If $f$ does not enjoy these properties on the whole interval  $[a,b]$, but does so on each subinterval~$[c_n,c_{n+1}]$ of a partition $c_0:=a<c_1<\cdots<c_N:=b$, we can write
\begin{IEEEeqnarray*}{rCl}
	I(\lambda) &=& \sum_{n=0}^{N-1}\int_{c_n}^{c_{n+1}}f(y)e^{i\lambda(y-\bar y)^2}dy,
\end{IEEEeqnarray*}
and handle each integral as before. In other words, only $a$, $b$, and the points where $f$ or its derivatives (below some prescribed order depending on the sought order of the expansion) is discontinuous or singular contribute to the expansion; notice that if we artificially introduce an intermediate point $c_n$ where $f$ is smooth enough, the contribution of this point will be zero since in this case $I_{c_n^+}(\lambda)+I_{c_n^-}(\lambda)=0$.

\begin{rem}\label{rem:asymp}
	There are marked differences between the cases $\bar y>b$ and $\bar y=b$. Indeed, assume the leading term in the asymptotic expansion of $f(b-r)$ is proportional to~$r^{\beta_0-1}$; when $\bar y>b$, the leading term of $\frac{(\chi_bf)(b-r)}{2(\bar b+r)}\bigr|_{r=-\bar b+\sqrt{\bar b^2+u}}$ is proportional to~$u^{\beta_0-1}$, whereas it is proportional to~$u^{\frac{\beta_0}{2}-1}$ when $\bar y=b$. Accordingly, the leading term of $I_b(\lambda)$ is proportional to~$e^{i\lambda(\bar y-b)^2}\lambda^{-\beta_0}$ when $\bar y>b$, and proportional to $\lambda^{-\frac{\beta_0}{2}}$ when $\bar y=b$. Similarly since $\bar y>a$, if the leading term in the asymptotic expansion of $f(a+s)$ is proportional to~$s^{\alpha_0-1}$, the leading term of $I_a(\lambda)$ is proportional to~$e^{i\lambda(\bar y-a)^2}\lambda^{-\alpha_0}$.
\end{rem}
\begin{rem}
	There is a similar but more complicated result, see~\cite[section~6.3]{BleistH1986book}, when, instead of ~\eqref{eq:Fasymptotic}, $F$ has the more general asymptotic expansion
	\begin{IEEEeqnarray*}{rCl}
		F(u) &=& \sum_{m=0}^M\sum_{n=0}^{N(m)}p_{mn}u^{a_m}(\log u)^n+o\bigl(u^{a_M}(\log u)^M\bigr).
	\end{IEEEeqnarray*}
\end{rem}

We can thus compute~\eqref{eq:uFirst} by setting $f(y):=\frac{\theta_0^{odd}(y)}{\sqrt{4\pi it}}$ and
$\lambda:=\frac{1}{4t}$; $\theta_0^{odd}$ is assumed smooth enough on each subinterval $[c_n,c_{n+1}]$ of a partition $c_0:=-1<c_1<\cdots<c_N:=1$, with well-behaved asymptotic expansions as $y\to c_n^+$ and $y\to c_{n+1}^-$. Notice that since by assumption $\theta_0\in L^2(0,1)$, the leading term in the asymptotic expansion of $\theta_0^{odd}(c_n+s)$ (resp. of $\theta_0^{odd}(c_n-r)$) is proportional to~$s^{\alpha_{0n}-1}$, with $\alpha_{0n}>\frac{1}{2}$ (resp. proportional  to~$r^{\beta_{0n}-1}$, with $\beta_{0n}>\frac{1}{2}$). In view of Remark~\ref{rem:asymp}, this means the leading term in the asymptotic expansion of~\eqref{eq:uFirst} due to~$c_n^+$, $i=0,\ldots,N-1$, is proportional to~$e^{i\frac{(1-c_n)^2}{4t}}t^{\alpha_{0n}-\frac{1}{2}}$ (resp. due to~$c_n^-$, $i=1,\ldots,N-1$ is proportional to~$e^{i\frac{(1-c_n)^2}{4t}}t^{\beta_{0n}-\frac{1}{2}}$), hence tends to~0 in an oscillatory way. On the other hand, the leading term  due to~$c_N=1$ is proportional to~$t^{\frac{\beta_{0N}-1}{2}}$, hence is not oscillatory but possibly singular if $\theta_0$ is itself singular at~$y=1$; since $\frac{\beta_{0N}-1}{2}>-\frac{1}{4}$, this term nonetheless belongs to~$L^4(\Rset)$, which is of course in accordance with \eqref{it:e} of Proposition~\ref{prop2}.

\subsection{Second phase}\label{sec:Sphase}
To compute sufficiently many terms in the series~\eqref{eq:control2}, an efficient way is to proceed recursively; moreover, the various terms must be properly scaled to accommodate finite precision arithmetics.

Following appendix~\ref{app:Gevrey}, we first compute
\begin{IEEEeqnarray*}{rClCl}
	\tilde\phi_k(t) &:=& \frac{(T-\tau)^k}{(2k)!}\frac{d^k}{dt^k}\Bigl[\phi_s\Bigl(\frac{t-\tau}{T-\tau}\Bigr)\Bigr]\\
	&=& \frac{1}{(2k)!}\phi_s^{(k)}\Bigl(\frac{t-\tau}{T-\tau}\Bigr).
\end{IEEEeqnarray*}
We emphasize that $\phi_s^{(k)}\bigl(\frac{t-\tau}{T-\tau}\bigr)$ is the $k^{th}$ derivative of $\phi_s$ evaluated at $\frac{t-\tau}{T-\tau}$, whereas $\frac{d^k}{dt^k}\bigl[\phi_s\bigl(\frac{t-\tau}{T-\tau}\bigr)\bigr]$ is the $k^{th}$ derivative of $t\mapsto\varphi_s(t):=\phi_s\bigl(\frac{t-\tau}{T-\tau}\bigr)$, i.e., $\varphi_s^{(k)}(t)$.

We next compute
\begin{IEEEeqnarray*}{rCl}
	\tilde y_k(t) &:=& \frac{\tau^k}{k!}\partial_t^k\theta^-_x(t,0).
\end{IEEEeqnarray*}
It can be checked directly that $\partial_xE$ satisfies the recurrence relation
\begin{IEEEeqnarray*}{rCl}
	\partial_t^{k+1}\partial_xE(t,x)&=& 
	-\frac{1}{4t^2}\bigg\{\bigl(ix^2+2t(4k+3)\bigr)\partial_t^{k}\partial_xE(t,x)\\
	&&\qquad+\>2k(2k+1)\partial_t^{k-1}\partial_xE(t,x)\bigg\};
\end{IEEEeqnarray*}
this can also be seen as a consequence of \eqref{eq:Hgenerating} and~\eqref{eq:Hrecurrence}.
Since $\partial_xE(t,x)=\frac{ix}{2}F(t,x)$, where $F(t,x):=\frac{E(t,x)}{t}$, this immediately implies
\begin{IEEEeqnarray}{rCl}
	\tilde F_{k+1}(t,x) &=& 
	\frac{-\tau}{4(k+1)t^2}\Bigl\{\bigl(ix^2+2t(4k+3)\bigr)\tilde F_k(t,x)\nonumber\\
	&&\qquad+\>2(2k+1)\tau\tilde F_{k-1}(t,x)\Bigr\},\label{eq:Ftilde}
\end{IEEEeqnarray}
where $\tilde F_k(t,x):=\frac{\tau^k}{k!}\partial_t^{k}F(t,x)$. On the other hand,
\begin{IEEEeqnarray*}{rCl}
	\theta^-_x(t,0) &=& -\int_{-1}^1\partial_yE(t,-y)\theta_0^{odd}(y)dy\\
	&=& -\frac{i}{2}\int_{-1}^1\frac{E(t,y)}{t}y\theta_0^{odd}(y)dy\\
	&=& -i\int_0^1F(t,y)y\theta_0(y)dy,
\end{IEEEeqnarray*}
since the integrand in the second line is even. Finally,
\begin{IEEEeqnarray*}{rCl}
	\tilde y_k(t) &=& -i\int_0^1\tilde F_k(t,y)y\theta_0(y)dy,
\end{IEEEeqnarray*}
with $\tilde F_k$ given by~\eqref{eq:Ftilde}. As $t>\tau$, performing the numerical integration causes no practical problem, unless $\tau$ is chosen very small.

We finally compute the scaled derivatives of~$Y$ using the general Leibniz rule,
\begin{IEEEeqnarray*}{rCl}
	\frac{Y^{(j)}(t)}{(2j)!} &=& 
	\frac{1}{(2j)!}\sum_{k=0}^j\binom{j}{k}\partial_t^k\theta^-_x(t,0)\frac{d^{j-k}\phi_s}{dt^{j-k}}\Bigl(\frac{t-\tau}{T-\tau}\Bigr)\\
	&=& \sum_{k=0}^j\frac{(2j-2k)!j!}{(2j)!(j-k)!\tau^k}\tilde y_k(t)\tilde\phi_k(t);
\end{IEEEeqnarray*}
the coefficient $d^j_k:=\frac{(2j-2k)!j!}{(2j)!(j-k)!\tau^k}$ is best computed recursively by $d^j_0=1$ and $d^j_k=\frac{d^j_{k-1}}{2\tau(2j-2k+1)}$.

The control in the second phase is eventually obtained by truncating the series~\eqref{eq:control2} and expressing it in terms of the scaled derivatives of~$Y$, which yields
	\begin{IEEEeqnarray}{rCl}
	u(t) &=& \sum_{j=0}^{\bar{\jmath}}\frac{(-i)^j}{(2j+1)}\frac{Y^{(j)}(t)}{(2j)!}.\label{eq:truncated}
\end{IEEEeqnarray}

\begin{figure}
	\begin{center}
		\includegraphics[width=\columnwidth]{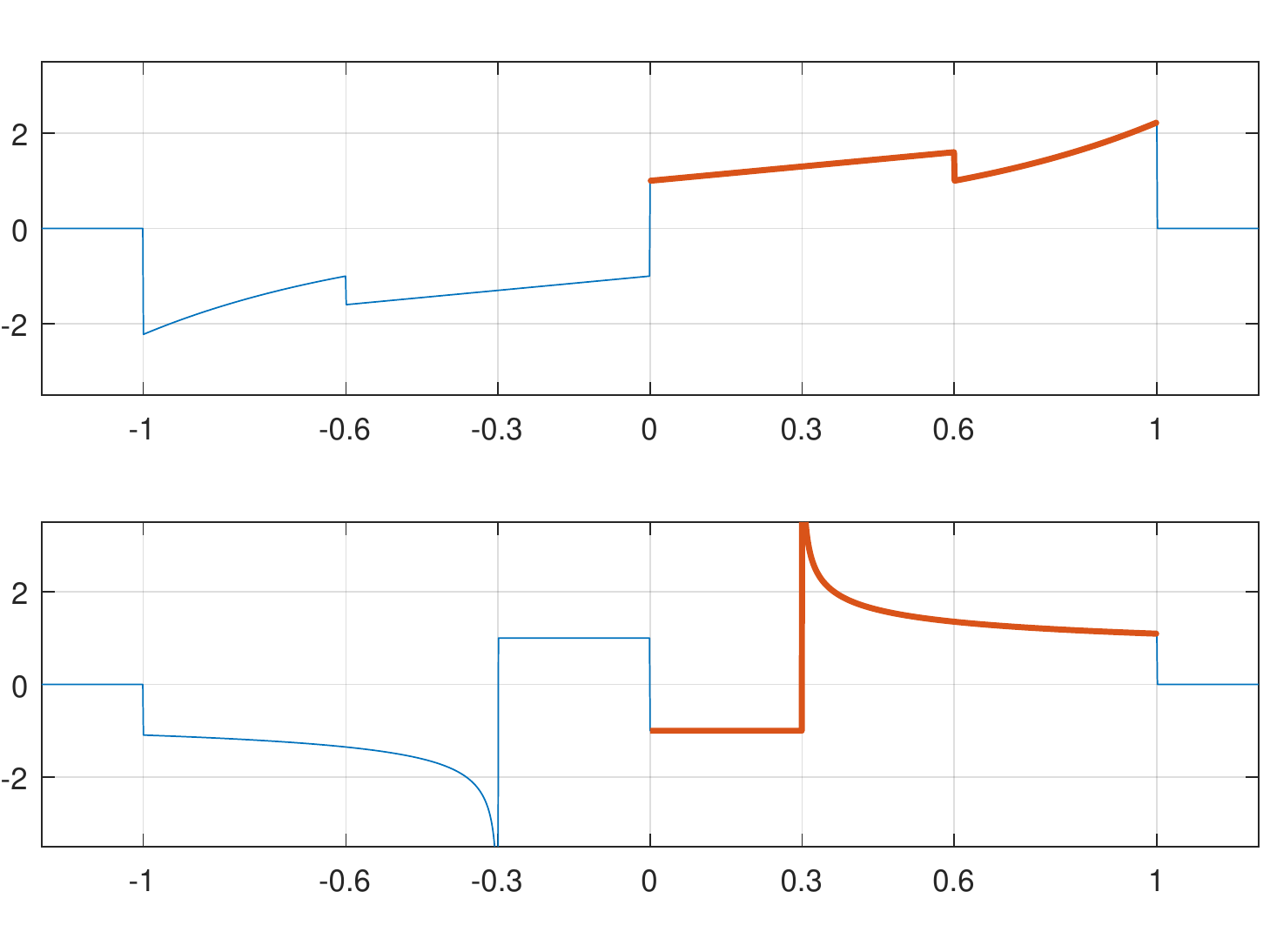}
		\caption{Initial condition~$\theta_0$ (red) and odd extension~$\theta_0^{odd}$ (blue); real parts (top), imaginary parts (bottom).}
		\label{fig:IC}
	\end{center}
\end{figure}

\section{Numerical experiments}\label{sec:simulation}
We now illustrate the effectiveness of the approach on a numerical example. The Matlab source code is freely available  as ancillary files to the preprint~\cite{MartiRR2017arXiv170501052}.
We choose as initial condition
\begin{IEEEeqnarray*}{rCl}
	\theta_0(x) &:=& \begin{dcases*}
		x+1-i, & $0<x\le0.3$\\
		x+1+i(x-0.3)^{-1/4}, & $0.3<x\le0.6$\\
		e^{2(x-0.6)}+i(x-0.3)^{-1/4}, & $0.6<x\le1$,
	\end{dcases*}
\end{IEEEeqnarray*}
displayed in Fig.~\ref{fig:IC}, and we want to steer the system to zero at time~$T:=0.4$, with intermediate time~$\tau:=0.05$. Notice the initial condition is rather challenging, since it does not satisfy the boundary conditions, has discontinuities, and a singularity (it can be shown to belong only to~$H^s(0,1)$, $s<\tfrac{1}{4}$).

\begin{figure}
	\begin{center}
		\includegraphics[width=\columnwidth]{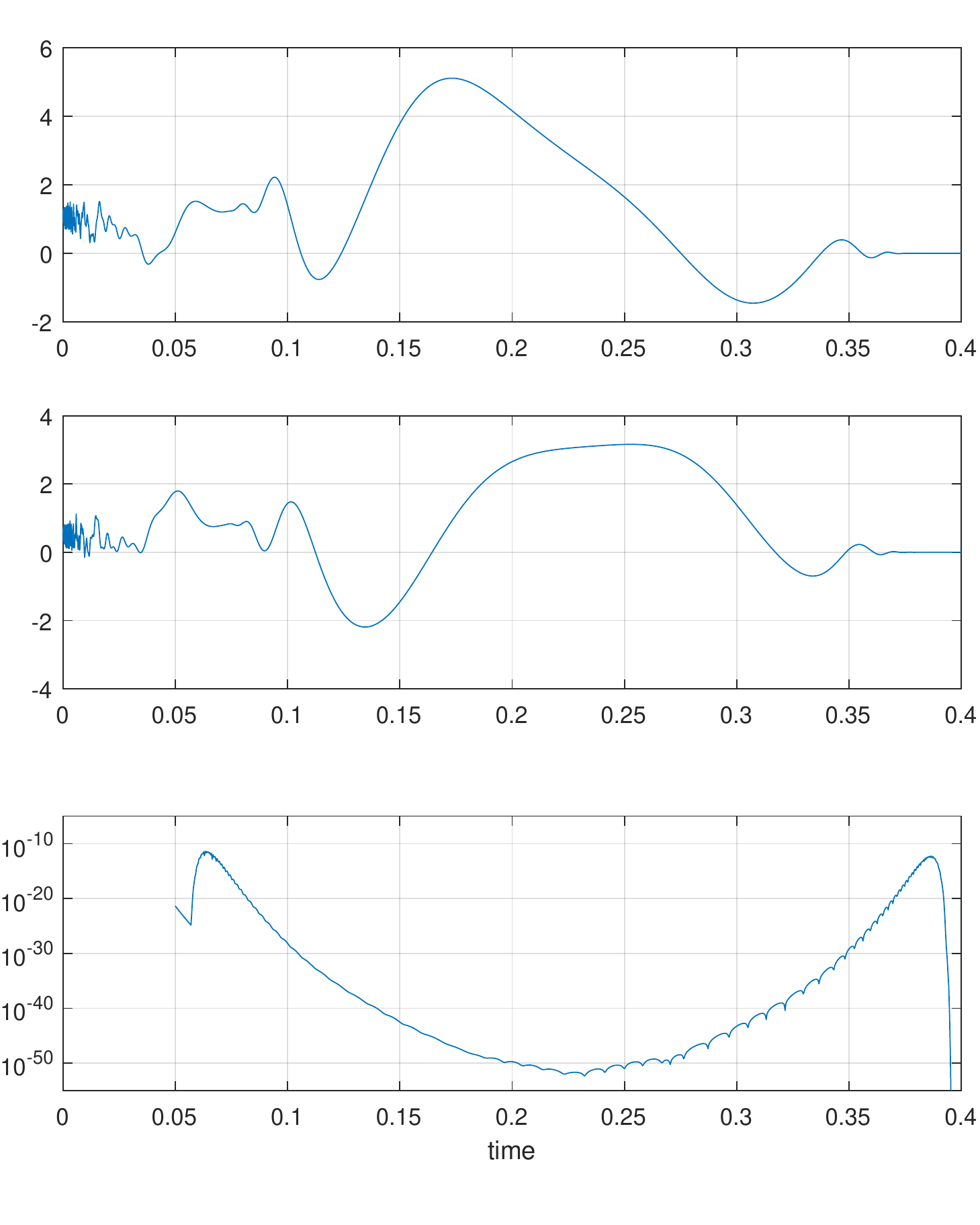}
		\caption{Control~$u(t)$: real part (top) and imaginary part (middle); contribution of 30 following terms in truncated control~\eqref{eq:truncated} (bottom).}
		\label{fig:control}
	\end{center}
\end{figure}

\begin{figure}
	\begin{center}
		\includegraphics[width=\columnwidth]{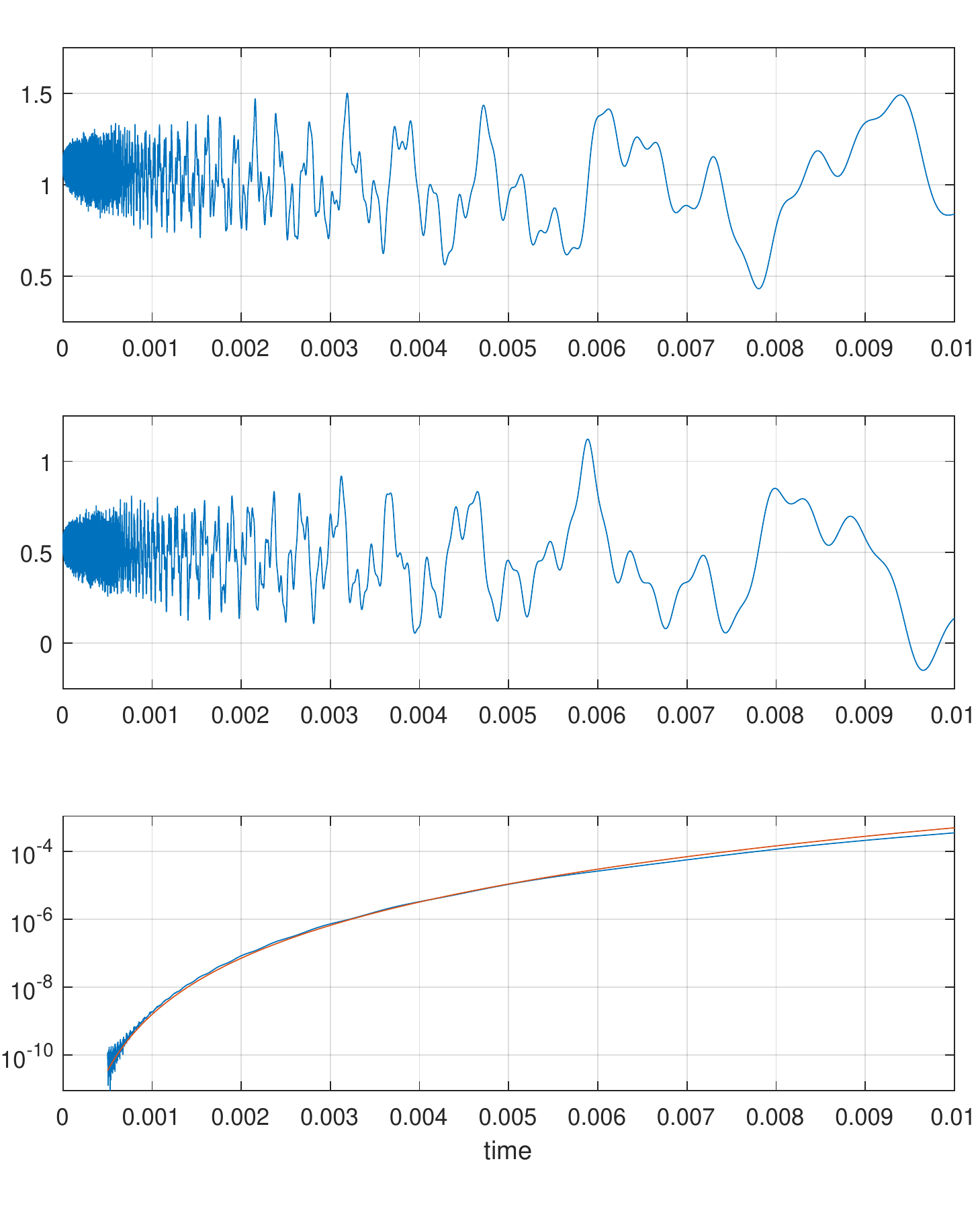}
		\caption{Control~$u(t)$ at the beginning of first phase: real part (top) and imaginary part (middle). Bottom: error between asymptotic expansion and numerical integration (blue), graph of $t\mapsto\num{5e7}t^\frac{11}{2}$ (red).}
		\label{fig:controlZoom}
	\end{center}
\end{figure}

The control, computed as explained in section~\ref{sec:implementation}, is shown in Fig.~\ref{fig:control}.
For the second phase, the control is given by~\eqref{eq:truncated}, with $s:=1.7$, $M=0.8$ and $\bar{\jmath}:=50$; the figure also shows that the contribution of the 30 following terms in the series is very small. For the first phase, the integral~\eqref{eq:uFirst} is evaluated by an asymptotic expansion for $0\le t<\num{e-3}$, and by a numerical quadrature (with Matlab \texttt{quadgk}) for $\num{e-3}\le t<\tau$; the expansion uses sufficiently many terms so that the neglected terms have order~$\tfrac{11}{2}$. 
Fig.~\ref{fig:controlZoom} displays the beginning of the first phase, together with the error between the ``exact'' values computed by numerical integration and the values obtained by expansion. The error is as expected of order~$\tfrac{11}{2}$; notice \texttt{quadgk} begins to have trouble evaluating the integral with a good accuracy when $t<\num{e-3}$. The control in the first phase, which is very oscillatory at the beginning because of the discontinuities and singularities of~$\theta_0^{odd}$, gradually calms down.

\begin{figure}
	\begin{center}
		\includegraphics[width=\columnwidth]{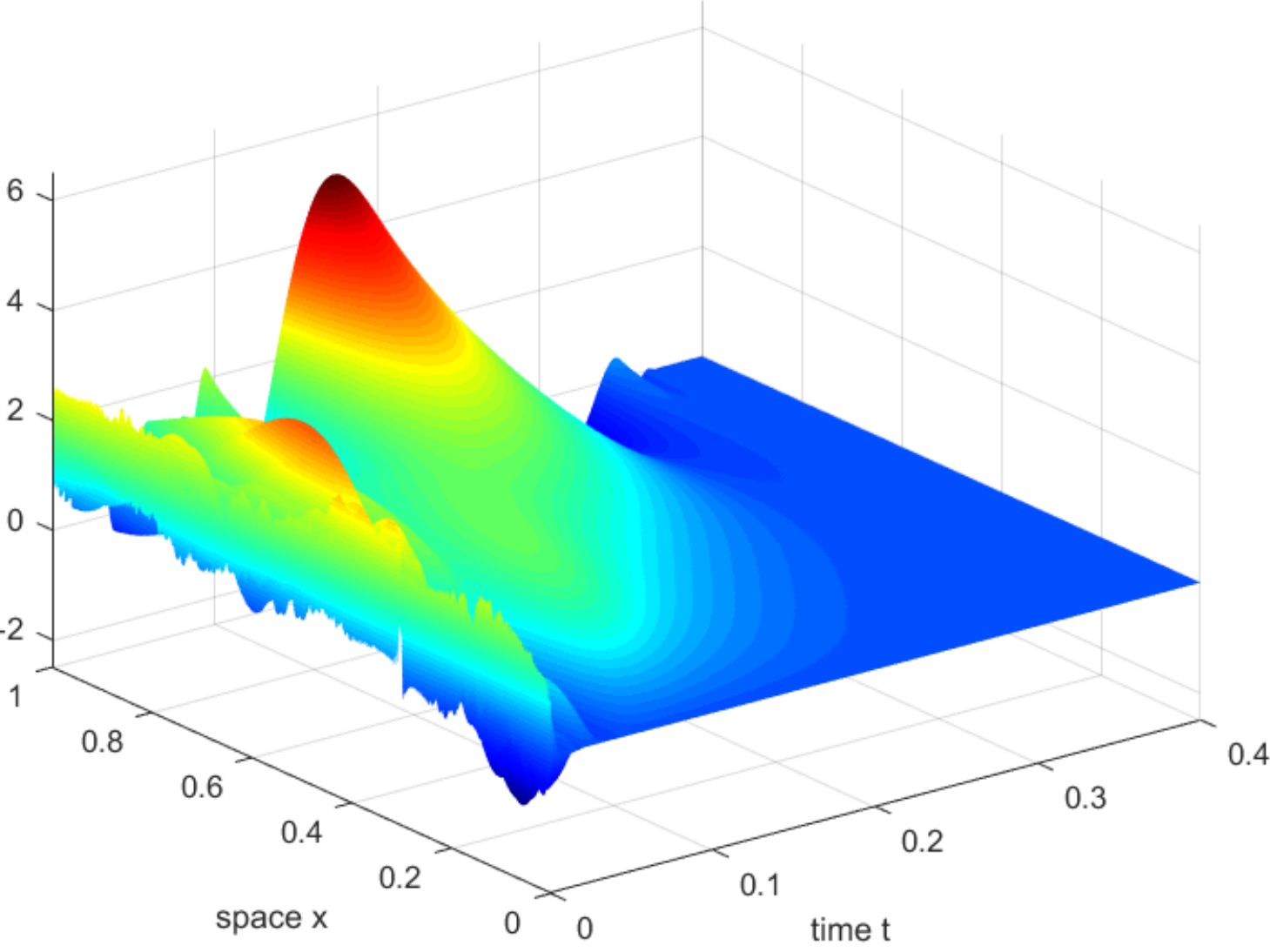}
		\includegraphics[width=\columnwidth]{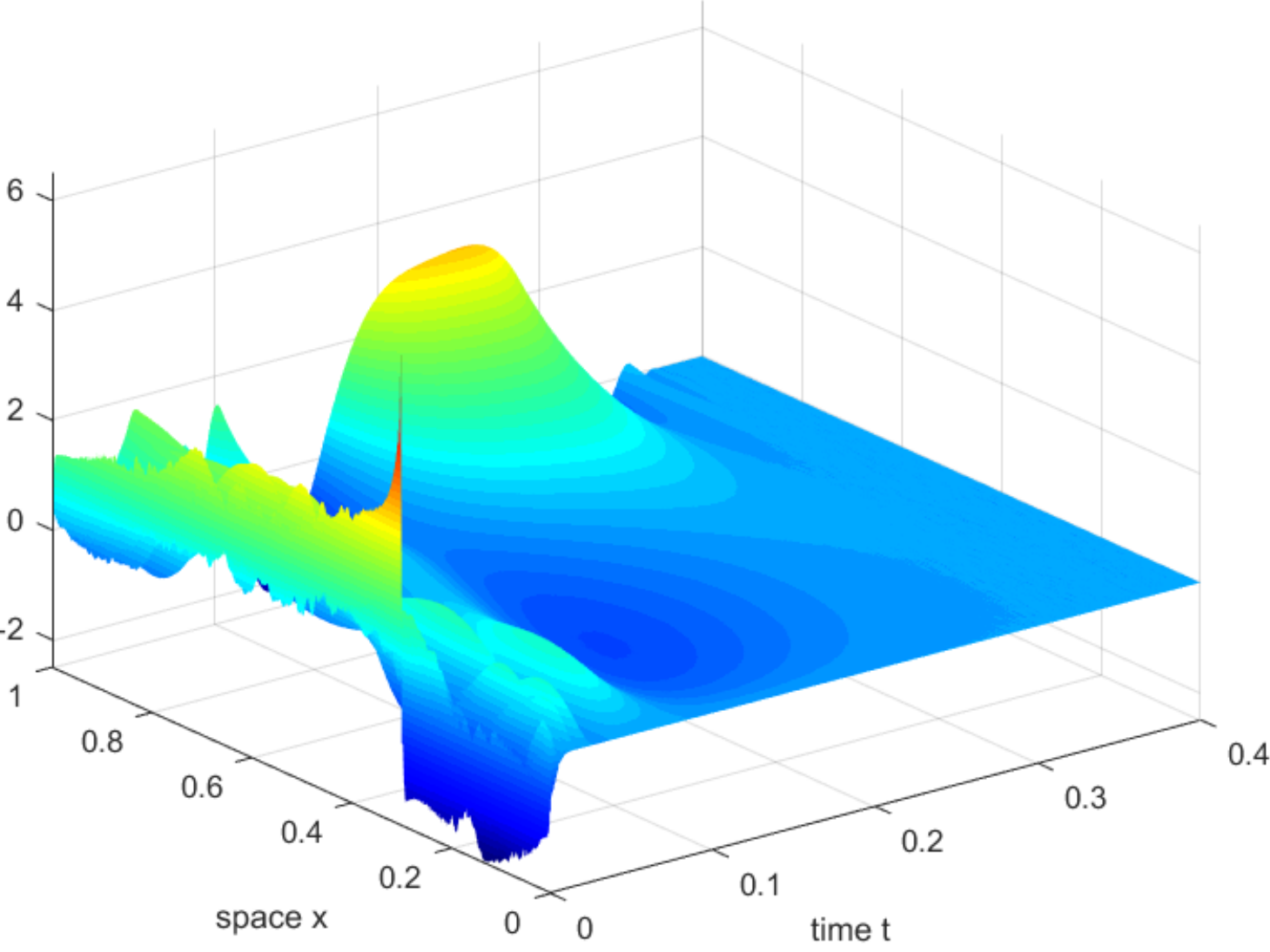}
	\end{center}
		\caption{Evolution of~$\theta$: real part (top) and imaginary part (bottom).}
\label{fig:3D}
\end{figure}

\begin{figure}
	\begin{center}
		\includegraphics[width=\columnwidth]{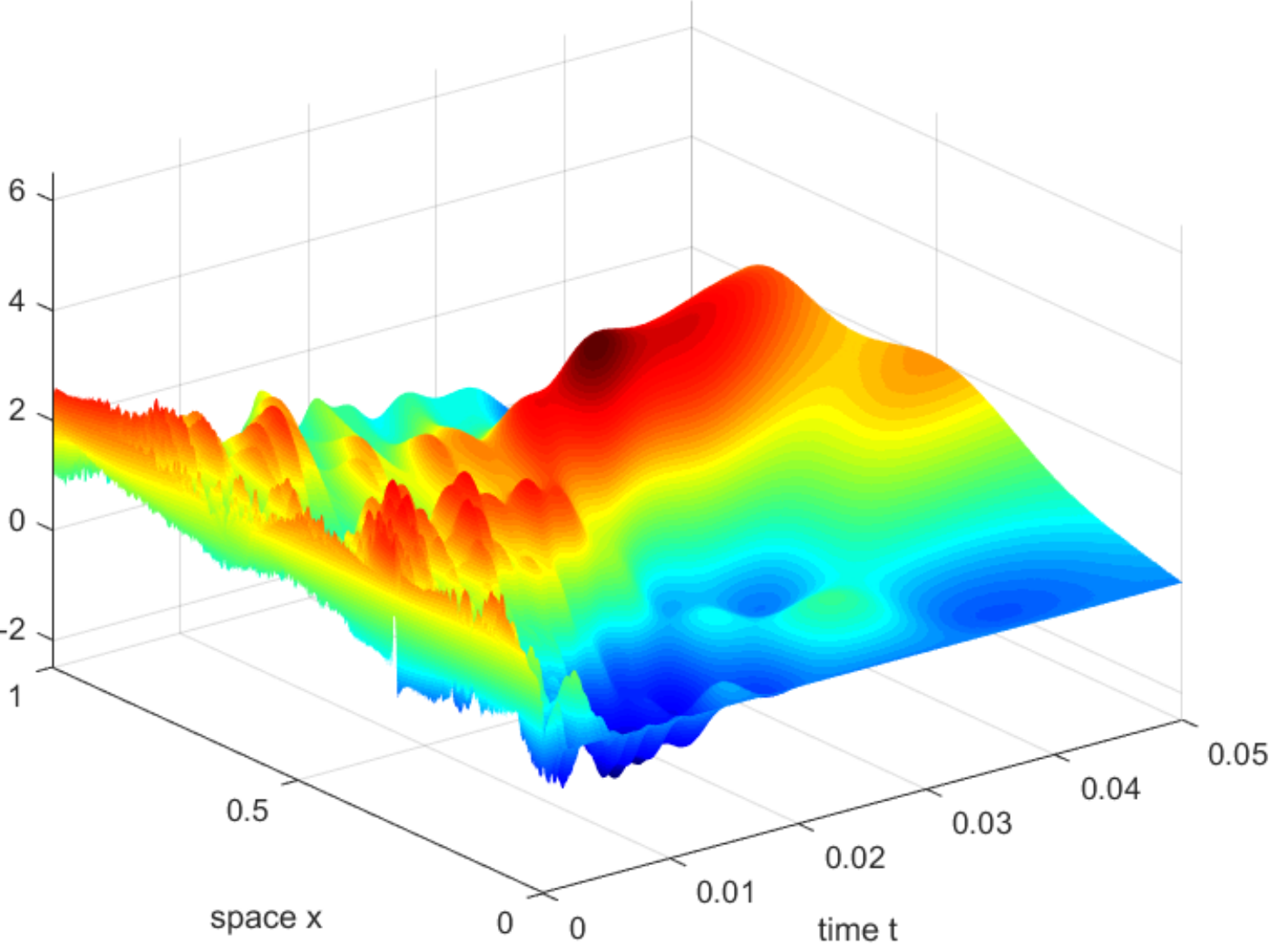}
		\includegraphics[width=\columnwidth]{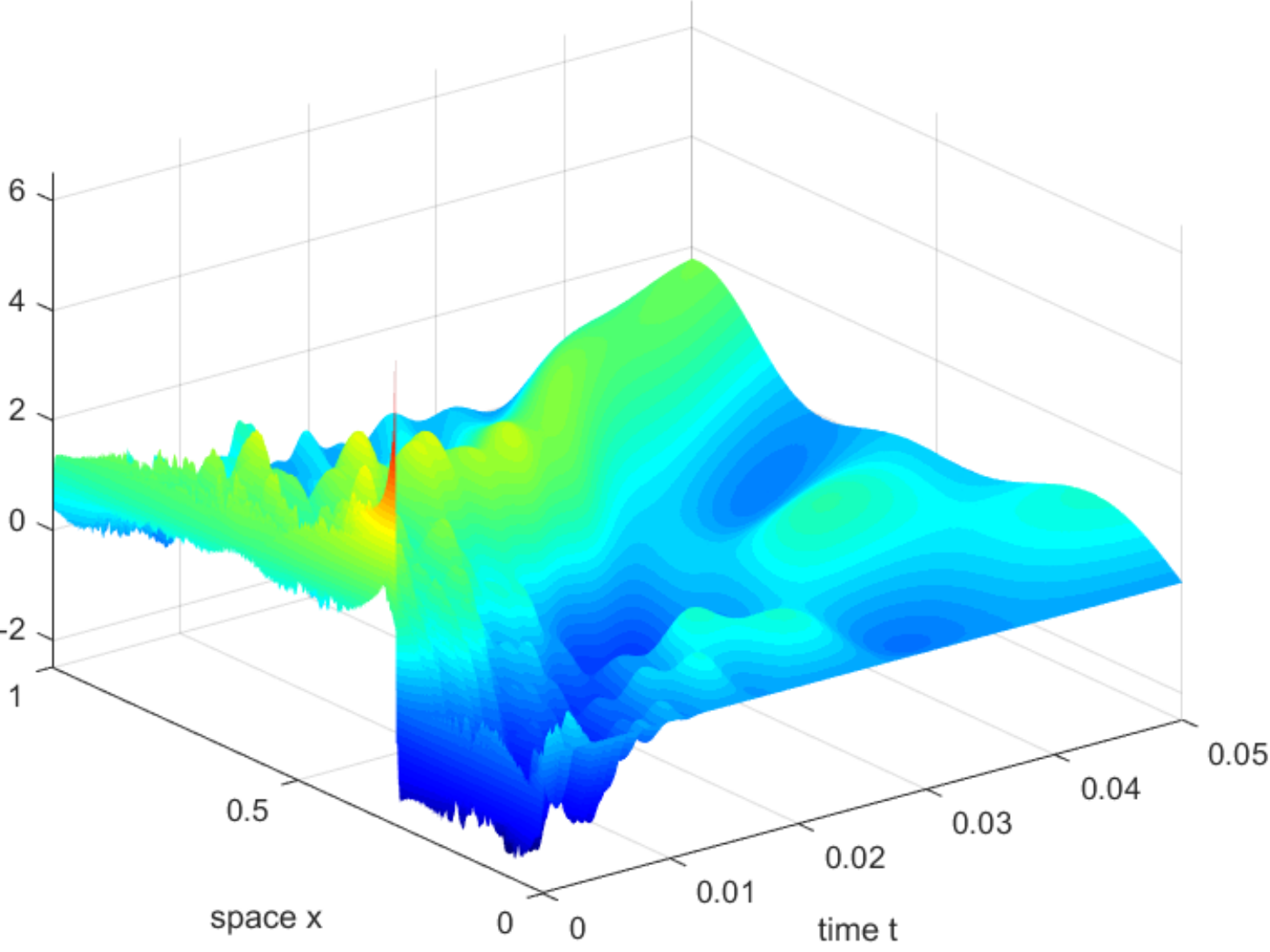}
	\end{center}
	\caption{Evolution of~$\theta$, zoom on first phase: real part (top) and imaginary part (bottom).}
	\label{fig:3DZoomFirst}
\end{figure}

\begin{figure}
	\begin{center}
		\includegraphics[width=\columnwidth]{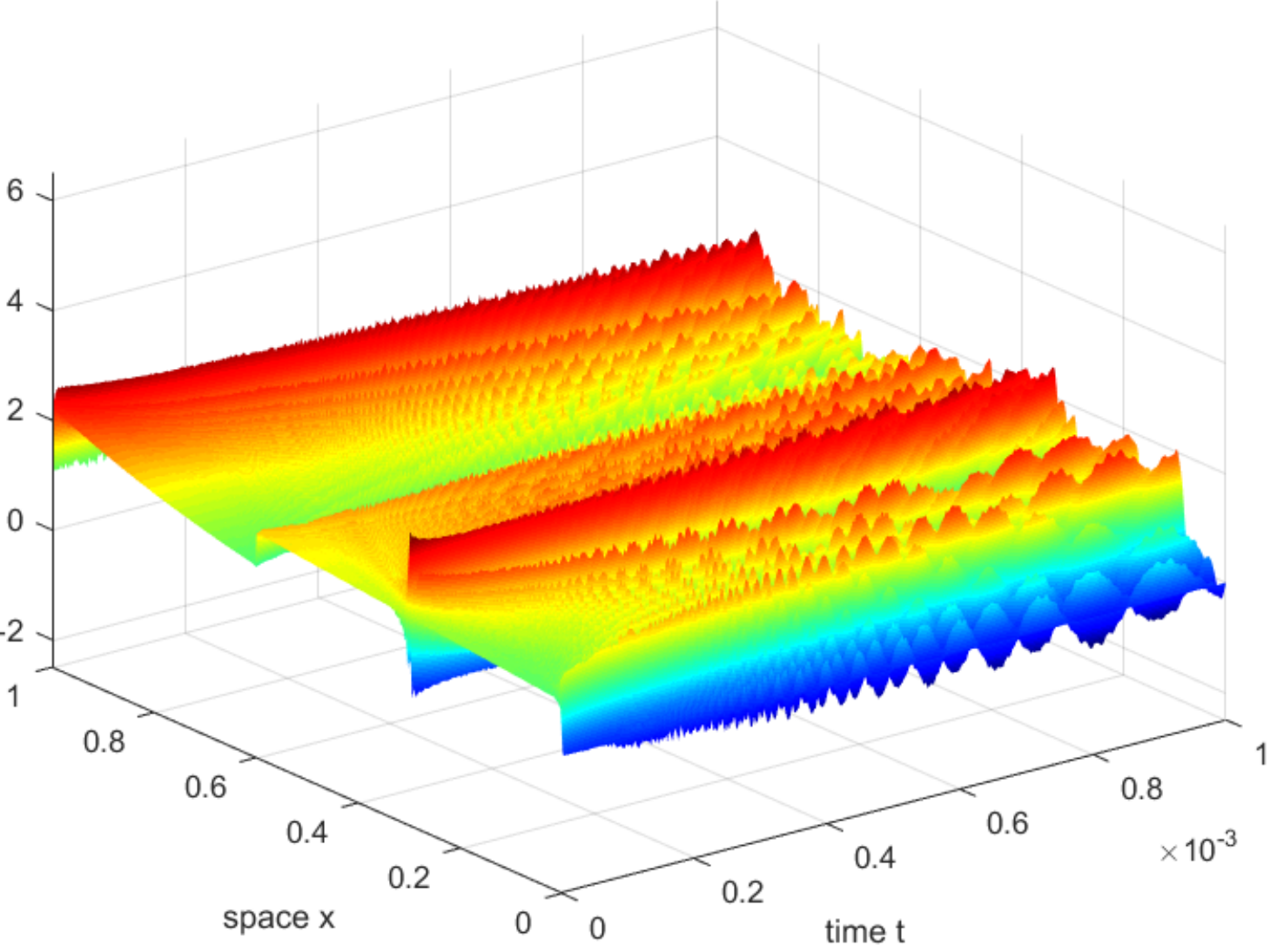}
		\includegraphics[width=\columnwidth]{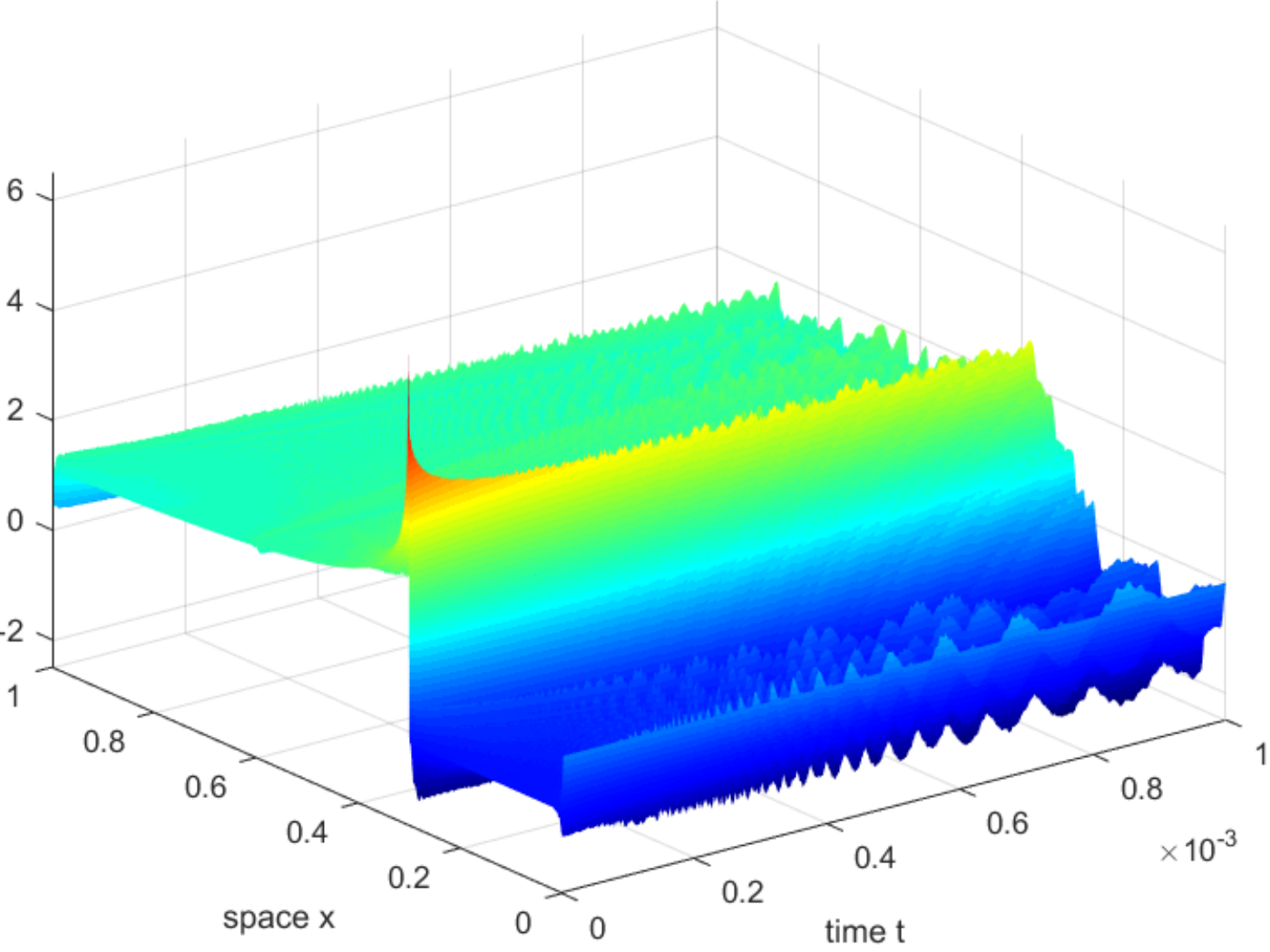}
	\end{center}
	\caption{Evolution of~$\theta$, zoom on beginning of first phase: real part (top) and imaginary part (bottom).}
	\label{fig:3DZoomInit}
\end{figure}

This control is then applied to the numerically simulated system. The numerical scheme is a fixed-step Crank-Nicolson scheme; due to the oscillatory nature of the Schrödinger equation, a small diffusion term $dx^\frac{3}{4}$, where $dx$ is the space step, is added to avoid spurious oscillations. For the same reasons, it is necessary to use very fine space and time grids to get a good accuracy. The scheme is initialized with four backward Euler half-steps to better handle the discontinuities in the initial condition (so-called ``Rannacher time-stepping''~\cite{Ranna1984NumMat}).
Fig.~\ref{fig:3D} displays the whole evolution of the system; the final error is about~\num{1.5e-3}. Fig.~\ref{fig:3DZoomFirst} zooms on the first phase; it can be seen that the control has as anticipated a smoothing effect, after a very oscillatory transient. Fig.~\ref{fig:3DZoomInit} zooms further on the very beginning of the motion; it shows how the discontinuities and singularity of the initial condition generate the initial oscillatory behavior.

The influence of~$\tau$ on the shape of the trajectory, as well as the influence of the order and type of the Gevrey step, is an interesting and open question; see~\cite{MartiRR2015ECC} for some numerical experiments in that regard.


\bibliographystyle{abbrv}        
\bibliography{schro1D}        

\appendix
\section{Recursive computation of $\frac{r^j\phi_s^{(j)}}{(2j)!}$}\label{app:Gevrey}
A practical problem when implementing the control in the second phase is to evaluate sufficiently many derivatives of~$\phi_s$, see section~\ref{sec:Sphase}. An efficient way is to proceed recursively, with a suitable scaling to accommodate finite precision arithmetics.

Let $u(\rho):=-\frac{M}{\rho^\sigma}$. Then $\rho u'=-\sigma u$, and applying the general Leibniz rule yields
\begin{IEEEeqnarray*}{rCl}
	\rho u^{(j)} + (j-1)u^{(j-1)} &=& -\sigma u^{(j-1)}.
\end{IEEEeqnarray*}
After scaling, this gives at once
\begin{IEEEeqnarray*}{rCl}
	\rho\frac{r^ju^{(j)}}{j!} &=& -r\biggl(1+\frac{\sigma-1}{j}\biggr) \frac{r^{j-1}u^{(j-1)}}{(j-1)!}.
\end{IEEEeqnarray*}
Similarly, the derivatives of $v(\rho):=-\frac{M}{(1-\rho)^\sigma}$ satisfy
\begin{IEEEeqnarray*}{rCl}
	(1-\rho)\frac{r^jv^{(j)}}{j!} &=& r\biggl(1+\frac{\sigma-1}{j}\biggr) \frac{r^{j-1}v^{(j-1)}}{(j-1)!}.
\end{IEEEeqnarray*}

Consider now $f(\rho):=e^{u(\rho)}$. Then $f'=u'f$, and applying the general Leibniz rule yields
\begin{IEEEeqnarray*}{rCl}
	f^{(j)} &=& \sum_{k=0}^{j-1}\binom{j-1}{k}u^{(k+1)}f^{(j-1-k)}\\
	&=& \sum_{k=1}^{j}\binom{j-1}{k-1}u^{(k)}f^{(j-k)}.
\end{IEEEeqnarray*}
After scaling, this gives at once
\begin{IEEEeqnarray*}{rCl}
	\frac{r^jf^{(j)}}{(2j)!} &=& \sum_{k=1}^{j}c^j_k \frac{r^ku^{(k)}}{k!} \frac{r^{j-k}f^{(j-k)}}{(2j-2k)!},
\end{IEEEeqnarray*}
where $c^j_k := k \frac{(2j-2k)!}{(2j)!} \frac{(j-1)!}{(j-k)!}$. The $c^j_k$ can be computed recursively by $c^j_1=\frac{1}{2j(2j+1)}$ and
$c^j_{k+1}=\frac{k+1}{2k(2j-2k-1)}c^j_k$.

Similarly, the scaled derivatives of $g(t):=e^{v(t)}$ read
\begin{IEEEeqnarray*}{rCl}
	\frac{r^jg^{(j)}}{(2j)!} &=& \sum_{k=1}^{j}c^j_k \frac{r^kv^{(k)}}{k!} \frac{r^{j-k}g^{(j-k)}}{(2j-2k)!}.
\end{IEEEeqnarray*}

Applying next the general Leibniz rule to $(f+g)\phi_s=g$ yields
\begin{IEEEeqnarray*}{rCl}
	(f+g)\phi_s^{(j)}+\sum_{k=1}^{j}\binom{j}{k}(f+g)^{(k)}\phi_s^{(j-k)} &=& g^{(j)}.
\end{IEEEeqnarray*}
After scaling,
\begin{IEEEeqnarray*}{rCl}
	(f+g)\frac{r^j\phi_s^{(j)}}{(2j)!} 
	&=& \frac{r^jg^{(j)}}{(2j)!} - \sum_{k=1}^{j}d^j_k\frac{r^kf^{(k)}+r^kg^{(k)}}{(2k)!}\frac{r^{j-k}\phi_s^{(j-k)}}{(2j-2k)!},
\end{IEEEeqnarray*}
where $d^j_k := \frac{(2k)!}{k!} \frac{(2j-2k)!}{(j-k)!} \frac{(j)!}{(2j)!}$.
The $d^j_k$ can be computed recursively by $d^j_0=1$ and $d^j_k=\frac{2k-1}{2j-2k+1}d^j_{k-1}$.

Setting $\tilde u^j:=\frac{r^ju^{(j)}}{j!}$, $\tilde v^j:=\frac{r^jv^{(j)}}{j!}$, $\tilde f^j:=\frac{r^jf^{(j)}}{(2j)!}$, and $\tilde g^j:=\frac{r^jg^{(j)}}{(2j)!}$, we finally obtain $\tilde\phi_s^j:=\frac{r^j\phi_s^{(j)}}{(2j)!}$ recursively by
\begin{IEEEeqnarray*}{rClrCl}
	\tilde u^0 &=& u,\qquad & \tilde u^j &=& -\frac{r}{\rho}\biggl(1+\frac{\sigma-1}{j}\biggr)\tilde u^{j-1}\\
	\tilde v^0 &=& v, & \tilde v^j &=& \frac{r}{1-\rho}\biggl(1+\frac{\sigma-1}{j}\biggr)\tilde v^{j-1}\\
	\tilde f^0 &=& f, & \tilde f^j &=& \sum_{k=1}^{j}c^j_k \tilde u^{k}\tilde f^{j-k}\\
	\tilde g^0 &=& g, & \tilde g^j &=& \sum_{k=1}^{j}c^j_k \tilde v^{k}\tilde g^{j-k}\\
	\tilde\phi_s^0 &=& \phi_s, & \tilde\phi_s^j &=& \frac{\tilde g^{j} - \sum_{k=1}^{j}d^j_k(\tilde f^{k}+\tilde g^{k})\tilde\phi_s^{j-k}}{\tilde f^0+\tilde g^0}.
\end{IEEEeqnarray*}
These formulas can then be directly implemented for instance in Matlab.


\end{document}